\documentclass[12pt]{amsart}
\usepackage{amsfonts}

\theoremstyle{plain}

\newtheorem{lemma}{Lemma}

\newtheorem{proposition}{Proposition}

\numberwithin{equation}{section}

\email{\footnotesize{nizar.demni@univ-rennes1.fr}}
\email{\footnotesize{mouayn@fstbm.ac.ma}}

\begin{document}
\author[N. Demni]{Nizar Demni}
\address{{\footnotesize IRMAR, Universit\'{e} de Rennes 1}\\
{\footnotesize \ Campus de Beaulieu}\\
{\footnotesize \ 35042 Rennes cedex}\\
{\footnotesize \ France}}
\author[Z. Mouayn]{Zouhair Mouayn}
\address{{\footnotesize Department of Mathematics}\\
{\footnotesize \ Faculty of Sciences and Technics (M'Ghila)}\\
{\footnotesize \ Sultan Moulay Slimane }\\
{\footnotesize \ PO. Box 523, B\'{e}ni Mellal}\\
{\footnotesize \ Morocco}}
\title{Analysis of generalized probability distributions associated with higher
Landau levels}
\maketitle

\begin{abstract}
To a higher Landau Level corresponds a generalization of the Poisson
distribution arising from generalized coherent states. In this paper, we
write down the atomic decomposition of this probability measure and
expressed its weights through ${}_2F_2$ hypergeometric polynomials. Then, we
prove that it is not infinitely divisible in opposite to the Poisson
distribution corresponding to the lowest Landau level. We also derive the
L\'evy-Kintchine representation of its characteristic function when the
latter does not vanish and deduce that the representative measure is signed.
By considering the total variation of the last measure, we obtain the
characteristic function of a new infinitely divisible discrete probability
distribution for which we compute also the
weights.
\end{abstract}


\section{Introduction}

\noindent The Poisson distribution is the cornerstone of the set of
integer-valued infinitely divisible distributions (\cite{Ste-VH}). In the
realm of quantum optics and precisely in photoemission, this probability
distribution describes the random number of photons per unit of time when
the intensity of the incident light is constant and the counts during
distinct time intervals are statistically independent (\cite{Kla-Sud},
p.16-19). It also arises from the \textit{standard coherent states} of the
Bargmann-Fock space of entire functions in the $L^{2}$-space of the complex
Gaussian distribution. In \cite{AIM}, this space was identified to the null
eigenspace of a magnetic Laplacian which is unitarily equivalent to the
Landau operator. The latter represents the Hamiltonian of a charged particle
evolving in the plane subject to the action of a homogeneous normal magnetic
field. More generally, the spectrum of the Landau operator, known as
Euclidean Landau levels, was shown to be the set of nonnegative integers
and the corresponding eigenspaces were called generalized Bargmann-Fock
spaces (known also as true polyanalytic spaces, \cite{AF}). In \cite{Mou},
generalized coherent states were attached to each eigenspace and it was
proved in \cite{Mou-Tou} that for any $m\geq 0$, the number $X_{m}$ of
photons per unit of time associated with the $m$-th Landau level is governed
by the following probability distribution
\begin{equation}
\Pr \left( X_{m}=j\right) =\frac{\left( \min \left( m,j\right) \right) !}{%
\left( \max \left( m,j\right) \right) !}\lambda ^{\left| m-j\right|
}e^{-\lambda }\left( L_{\min \left( m,j\right) }^{\left( |m-j|\right)
}\left( \lambda \right) \right) ^{2},j\geq 0,  \label{PCD}
\end{equation}
where $L_{n}^{\left( \alpha \right) }$ is the $n$-th Laguerre polynomial of
parameter $\alpha >-1$ (\cite{Sze}) and $\lambda >0$ is the intensity of the
light. In particular, $X_{0}$ reduces to the Poisson distribution of
parameter $\lambda $. When $m\geq 1$, the characteristic function of $X_{m}$
was computed in (\cite{Mou-Tou}): it splits into the product of the
characteristic functions of $X_{0}$ and of a Dirac mass at $m$, together
with an additional oscillating factor expressed as a Laguerre polynomial
composed with a cosine function. Since this last factor is the Fourier
transform of a finitely-supported measure, it does not correspond to an
infinitely divisible distribution (\cite{Ste-VH}). However, this fact does
not imply that the distribution of $X_{m},m\geq 1$ is not infinitely
divisible and it is therefore natural to wonder whether this property
remains valid for higher Landau levels $m\geq 1$.

In this paper, we analyze the structure of the distribution of $X_m, m\geq 1$%
. More precisely, we write down its atomic decomposition and express its
weights through ${}_2F_2$ hypergeometric polynomials of $\lambda$. We also
prove that the infinite divisibility of $X_m$ breaks down as soon as $m\geq
1 $: even when the characteristic function of $X_m$ does not vanish, the
canonical sequence of the moment generating function of this random variable
has a signed first element (\cite{Ste-VH}). In order to get more insight
into the failure of this property, we derive a L\'evy-Kintchine
representation of the characteristic function of $X_m$ when it does not
vanish, whence it is readily seen that the representative measure is signed.
This fact not only proves again that $X_m$ is not infinitely divisible but
also provides a new infinitely divisible discrete probability distribution
whose L\'evy measure is the total variation of the last signed measure.
Using similar computations as t he previous ones leading to the
L\'evy-Kintchine representation, we compute the characteristic function of
the new distribution and write down its atomic decomposition.

The paper is organized as follows. In section 2, we recall the coherent
states formalism used to construct the generalized Poisson distribution
displayed in \eqref{PCD} and write down its atomic decomposition. In section
3, we prove the non infinite divisibility of $X_m$ and derive of the
aforementioned L\'evy-Kintchine-type representation. Section 4 is devoted to
the analysis of the newly obtained infinitely divisible distribution.

\section{Generalized Poisson distribution: atomic decomposition}

In this paragraph, we briefly recall from \cite{Mou-Tou} the construction of
the distribution of $X_{m}$ which we refer to below as the generalized
Poisson distribution. To this end, we introduce the Hamiltonian $H$ with
constant homogeneous magnetic field which may be written, up to
normalizations, as:
\begin{equation}
H:=-\frac{1}{4}\left( \left( \frac{\partial }{\partial x}+\sqrt{-1}y\right)
^{2}+\left( \frac{\partial }{\partial y}-\sqrt{-1}x\right) ^{2}\right) -%
\frac{1}{2}.
\end{equation}
Conjugating by the complex Gaussian function $\phi \mapsto e^{\frac{1}{2}%
|z|^{2}}\phi $, $\phi \in L^{2}\left( \mathbb{R}^{2},dxdy\right) $ where $%
z=x+\sqrt{-1}y$ then $H$ is unitarily mapped into the magnetic Laplacian
\begin{equation*}
L:=-\frac{\partial ^{2}}{\partial z\partial \overline{z}}+%
\overline{z}\frac{\partial }{\partial \overline{z}}.
\end{equation*}
This is a self-adjoint operator in the Hilbert space $L^{2}(\mathbb{C}%
,e^{-|z|^{2}}dz)$ where $dz$ is the Lebesgue measure in $\mathbb{C}$ and has
purely discrete positive spectrum given by the set of nonnegative integers (
\cite{AIM}). Actually, the direct sum decomposition holds
\begin{equation*}
L^{2}(\mathbb{C},e^{-|z|^{2}}dz)=\oplus _{m=0}^{\infty }A_{m}(%
\mathbb{C})
\end{equation*}
where
\begin{equation*}
A_{m}(\mathbb{C}):=\{Lf=mf\}\cap L^{2}(\mathbb{C}%
,e^{-|z|^{2}}dz)
\end{equation*}
is the eigenspace of $L$ corresponding to the $m$-th Landau level.
In particular, $A_{0}(\mathbb{C})$ coincides with the
Bargmann-Fock space of holomorphic functions in $L^{2}(\mathbb{C}%
,e^{-|z|^{2}}dz)$, for which an orthonormal basis is given by $(\psi
_{j}(z):=z^{j}/\sqrt{\pi j!},j\geq 0)$ yielding the reproducing kernel
\begin{equation*}
K_{0}(z,w)=\sum_{j\geq 0}\psi _{j}(z)\overline{\psi _{j}(w)}=\frac{1}{\pi }%
e^{\langle z,w\rangle },\quad z,w\in \mathbb{C}.
\end{equation*}
It follows that the sequence
\begin{equation*}
\frac{|\psi _{j}(z)|^{2}}{K_{0}(z,z)},\quad j\geq 0,
\end{equation*}
defines a discrete probability distribution which is nothing else but the
Poisson distribution of parameter $\lambda =|z|^{2}$. When $m\geq 1$, the
reproducing kernel of $A_{m}(\mathbb{C})$ was given closed formula
in \cite{AIM} and the basis elements coincide, up to a normalization, with
the complex Hermite polynomials (see \cite{Mou-Tou} where they were denoted $%
h_{m,j}$):
\begin{equation*}
\psi _{m,j}(z):=\frac{1}{\sqrt{\pi m!j!}}\sum_{n=0}^{\min (m,j)}(-1)^{n}%
\binom{m}{n}\frac{j!}{(j-n)!}z^{j-n}(\overline{z})^{m-n}.
\end{equation*}
Note in passing that these polynomials appears in relation with complex
multiple Wiener integrals (\cite{Ito}). Moreover, the generalized Poisson
distribution \eqref{PCD} was obtained in \cite{Mou-Tou} in the same fashion
the Poisson distribution was obtained from $A_{0}(\mathbb{C})$.
Independently and at nearly the same time, this distribution appeared in
\cite{CGG} where few statistical properties were derived. In \cite{Mou-Tou},
its characteristic functions was computed and takes the form:
\begin{equation}
\mathbb{E}[e^{iuX_{m}}]=\exp \left( \lambda \left( e^{iu}-1\right) \right)
e^{imu}L_{m}^{\left( 0\right) }\left( 2\lambda \left( 1-\cos u\right)
\right) ,\qquad u\in \mathbb{R}.  \label{CF}
\end{equation}
With this expression in hands, we perform the analysis of the generalized
Poisson distribution and write down in the subsequent paragraph its atomic
decomposition.

\subsection{Atomic decomposition}

Clearly, we only need to write $L_m^{(0)}[2\lambda(1-\cos u)]$ as a linear
combination of $\cos ju, 1 \leq j \leq m$. This is the content of the
following lemma.

\begin{lemma}
Let $m \geq 1$ and $\lambda >0$, then
\begin{equation}
L_{m}^{(0)}\left[2\lambda \left( 1-\cos u\right) \right] =
\gamma_0(\lambda,m) + 2\sum _{j=1}^{m} \gamma _{j}\left(\lambda, m\right)
\cos ju,  \label{6.24}
\end{equation}
where
\begin{eqnarray*}
\gamma _{j}\left(\lambda ,m\right) = (-1)^j\sum_{k=j}^m \binom{m}{k}\binom{2k%
}{k-j} \frac{(-\lambda)^k}{k!}.
\end{eqnarray*}
Moreover, $\gamma_j(\lambda,m)$ are ${}_2F_2$ hypergeometric polynomials of $%
\lambda$ (see \cite{Erd}, CH.IV for the definition of ${}_2F_2$).
\end{lemma}

\begin{proof}
The expression of $\gamma _{j}\left(\lambda ,m\right)$ follows from the expansion (\cite{Sze}, p.101)
\begin{equation*}
L_{m}^{(0)}\left[2\lambda \left(1-\cos u\right) \right] = \sum_{k=0}^m \binom{m}{k}\frac{(-2\lambda)^k}{k!} (1-\cos u)^k = \sum_{k=0}^m \binom{m}{k}\frac{(-\lambda)^k}{k!} 4^k\sin^{2k}(u/2)
\end{equation*}
together with the linearization formula (\cite{Gra-Ryz}, p.31)
\begin{equation}\label{Line}
4^k\sin ^{2k}(u/2)=\binom{2k}{k} + 2\sum_{j=1}^{k}\left( -1\right)^{j}\binom{2k}{k-j} \cos(ju).
\end{equation}
As to the last claim of the lemma, it suffices to change the index summation $k \mapsto m+k$ and to use the duplication formula (\cite{Erd}, p.5, eq.15):
\begin{align*}
(2j+2k)! = \frac{4^{j+k}}{\sqrt{\pi}}(j+k)!\Gamma(j+k+1/2)
 = 4^{j+k}(1/2)_j(j+1/2)_k(j+k)!
\end{align*}
where $(a)_k = a(a+1)\cdots(a+k-1)$ is the Pochhammer symbol.
More precisely
\begin{align*}
\gamma_j(\lambda,m) &= \frac{m!}{(m-j)!}\sum_{k=0}^{m-j}\frac{(-1)^k(m-j)!}{(m-j-k)!(k+j)!}\binom{2k+2j}{k}\frac{\lambda^{k+j}}{(k+j)!}
\\& = \frac{(1/2)_j(4\lambda)^j}{(2j)!}\binom{m}{j} \sum_{k=0}^{m-j}\frac{(j-m)_k}{(j+1)_k}\frac{(j+1/2)_k}{(2j+1)_k}\frac{(4\lambda^{k})}{k!}
\\&= \frac{\lambda^j}{j!}\binom{m}{j} {}_2F_2(j-m, j+1/2, j+1, 2j+1, 4\lambda).
\end{align*}
\end{proof}

From this lemma, it follows that:

\begin{proposition}
The distribution of $X_m$ is decomposed as
\begin{equation*}
\mathcal{P}\left( \lambda \right) \ast \delta _{m} \ast \left(
\sum\limits_{j=-m}^{m} \gamma _{|j|}\left( \lambda ,m\right) \delta_j\right)
= \mathcal{P}\left( \lambda \right) \ast \left( \sum\limits_{j=0}^{2m}
\gamma _{|m-j|}\left( \lambda ,m\right) \delta_j\right) := \sum_{j =
0}^{\infty}q_j(\lambda,m)\delta_j
\end{equation*}
where $\mathcal{P}\left( \lambda \right)$ stands for the Poisson
distribution of parameter $\lambda$ and
\begin{align*}
q_j(\lambda,m) := e^{-\lambda}\sum_{k=0}^{j \wedge 2m} \frac{\lambda^{j-k}}{%
(j-k)!}\gamma _{|m-k|}.
\end{align*}
\end{proposition}

\section{Infinite divisibility}

This section is devoted to the proof of the non infinite divisibility of $%
X_m $ as well as to the L\'evy-Kintchine representation for its
characteristic function when the latter does not vanish. Recall that a
random variable $Y$ is infinitely divisible if for any $n \geq 2$, there
exist $n$ independent and identically distributed random variables $%
(Y_i^{(n)})_{i=1}^n$ such that
\begin{equation*}
Y = \sum_{i=1}^nY_i^{(n)}.
\end{equation*}
A necessary condition for this property to hold is the non vanishing of the
characteristic function of the given random variable (\cite{Ste-VH}, CH.IV,
Proposition 2.4). Accordingly and with regard to \eqref{CF}, if $m\geq 1$ is
fixed then the parameter $\lambda $ must satisfy
\begin{equation}
2\lambda \left(1-\cos u\right) \neq x_{k}^{\left( m\right) },\qquad 1\leq
k\leq m,  \label{6.2}
\end{equation}
for all $u\in \mathbb{R}$, where $x_{k}^{\left( m\right) },1\leq k\leq m$
are the increasingly-ordered zeros of Laguerre polynomial $L_{m}^{\left(
0\right)}$, which are known to be simple and positive (\cite{Sze}). But
since
\begin{equation}
0\leq 2\lambda \left( 1-\cos u\right) \leq 4\lambda ,  \label{6.3}
\end{equation}
then the distribution of $X_m$ cannot be infinitely divisible unless
probably when $\lambda < x_{1}^{\left(m\right) }/4$. For those values of $%
\lambda$, this property does not hold neither:

\begin{proposition}
For any $m \geq 1$ and any $\lambda < x_{1}^{\left(m\right) }/4$, the
distribution of $X_m$ is not infinitely divisible.
\end{proposition}

\begin{proof}
Since $X_m$ is a integer-valued random variable, we can use for instance the following criterion (\cite{Ste-VH}, CH.II, Theorem 4.1) rather than other sophisticated criteria from the general theory: a discrete probability distribution $\left(p_{k}\right) _{k\in \mathbf{N}}$ with $p_{0}>0$ is infinitely divisible if and only if the logarithmic derivative of its moment generating function
\begin{equation}
P\left( z\right) :=\sum\limits_{k\geq 0}p_{k}z^{k},\quad z\in \left[ 0,1\right]   \label{6.5}
\end{equation}
satisfies
\begin{equation*}
\frac{P^{\prime }\left( z\right) }{P\left( z\right) }=\sum\limits_{k\geq 0}r_{k}z^{k}
\end{equation*}
for non negative reals $r_k, k \geq 0$. The sequence $(r_k)_{k \geq 0}$ is referred to as the canonical sequence of $P$ and we prove below that $r_1 < 0$. To proceed, we first deduce  the moment generating function of $X_m$ from its characteristic function \eqref{CF}. Actually, the latter extends to an analytic function in the upper-half plane $\mathbb{C}^{+}$ so that
\begin{equation*}
\mathbb{E}\left(e^{-uX_m}\right) =\exp \left( \lambda \left(e^{-u}-1\right) -mu\right) L_{m}^{(0)}\left( 2\lambda \left( 1-\cosh u\right)\right).
\end{equation*}
for any $u\geq 0$. Setting $z=e^{-u}$ then the moment generating function of $X_m$ reads
\begin{equation*}
P_m\left( z\right) := \mathbb{E}\left( z^{X_m}\right) =\exp \left(\lambda \left( z-1\right) \right) z^{m}L_{m}^{\left(0\right) }\left( -\frac{\lambda }{z}\left( 1-z\right) ^{2}\right).
\end{equation*}
Next, since the leading coefficient of $L_m^{(0)}$ is $(-1)^m$ then
\begin{align*}
\frac{P_m'\left( z\right) }{P_m\left( z\right) }& =\lambda +\frac{m}{z}+%
\frac{d}{dz}\ln L_{m}^{\left( 0\right) }\left( -\frac{\lambda }{z}\left(1-z\right) ^{2}\right)
\\& =\lambda +\frac{m}{z}+\frac{d}{dz}\sum_{k=1}^{m}\ln \left(x_{k}^{\left( m\right) }+\frac{\lambda }{z}\left( 1-z\right) ^{2}\right)
\\& =\lambda +\frac{m}{z}-\frac{1}{z}\sum_{k=1}^{m}\frac{\lambda \left(1-z^{2}\right) }{zx_{k}^{\left( m\right) }+\lambda \left( 1-z\right) ^{2}}
\\& =\lambda +\sum\limits_{k=1}^{m}\frac{x_{k}^{\left( m\right)}}{zx_{k}^{\left( m\right) }+\lambda \left(1-z\right)^{2}}
\\& = \lambda +\sum\limits_{k=1}^{m}\frac{x_{k}^{\left( m\right)}}{1+(x_k^{(m)} - 2\lambda)z+z^2}.
\end{align*}
By assumption $x_k^{(m)} - 2\lambda > 0$ therefore
\begin{equation*}
r_1 = \frac{d}{dz} \left(\frac{P_m'}{P_m}\right)(0) = \sum_{k=1}^mx_k^{(m)}(2\lambda-x_k^{(m)}) < 0.
\end{equation*}
\end{proof}

\subsection{A L\'{e}vy-Khintchine representation}

In order to get more insight into the failure of the non infinite
divisibility property of $X_m$ under the assumption $4\lambda < x_1^{(m)}$,
we shall derive a L\'evy-Kintchine representation for its characteristic
function. Recall from \cite{Sato} (see p.37, Theorem 8.1) that if $Y$ is an
infinitely divisible random variable without Gaussian component, then its
L\'evy-Kintchine representation may be written as
\begin{equation*}
\ln \mathbb{E}[e^{iuY}] = ibu + \int (e^{iux} - 1 - iux\mathbf{1}_{|x| < 1})
\nu(dx)
\end{equation*}
where $b \in \mathbb{R}$ and $\nu$, the L\'evy measure of $Y$, is a positive
measure satisfying
\begin{equation}  \label{LM}
\nu\{0\} = 0, \quad \quad \int (1 \wedge x^2) \nu(dx) \, < \, \infty.
\end{equation}
Coming back to $X_m$, it is obvious that
\begin{equation*}
\lambda \left( e^{iu}-1\right) +imu = imu + \lambda \int (e^{iux} - 1)
\delta_1(dx),
\end{equation*}
so that we are only concerned with the representation of
\begin{align*}
\ln \left(L_{m}^{\left( 0\right) }\left( 2\lambda \left(1-\cos u\right)
\right) \right) &= \sum_{k=1}^{m}\ln \left(x_{k}^{\left(m\right)} - 4\lambda
\sin ^{2}(u/2)\right) = \sum_{k=1}^m \ln x_{k}^{(m)} +
\sum_{k=1}^m\ln\left(1-a_{k}^{(m)}\sin^{2}(u/2)\right)
\end{align*}
where
\begin{equation*}
a_{k}^{(m)} :=\frac{4\lambda }{x_{k}^{(m)}} \in (0,1), \, 1 \leq k \leq m.
\end{equation*}
In this respect, we prove

\begin{proposition}
Let $m \geq 1$. Then, for every $1 \leq k \leq m$, there exists a signed
measure $\mu_k^{(m)}$ satisfying \eqref{LM} and such that
\begin{equation*}
\ln \left(L_{m}^{\left( 0\right) }\left( 2\lambda \left(1-\cos u\right)
\right) \right) = \int(e^{iux}-1) \left[\sum_{k=1}^m \mu_k^{(m)}(dx)\right].
\end{equation*}
\end{proposition}

\begin{proof}
Since $L_m^{(0)}(x) = (-1)^mx^m + $ terms of lower degrees and  $L_m^{(0)}(0) = 1$ then $\prod_{i=1}^{m}x_{i}^{(m)}=1$ whence
\begin{equation*}
\ln \left(L_{m}^{\left( 0\right) }\left( 2\lambda \left(1-\cos u\right) \right) \right) = \sum_{k=1}^m \ln\left(1-a_{k}^{(m)}\sin^{2}(u/2)\right).
\end{equation*}
Now, fix $k \in \{1,\dots, m\}$ and expand
\begin{align}\label{E1}
\ln \left( 1-a_{k}^{(m)}\sin ^{2}\frac{u}{2}\right) &= - \sum_{j=1}^{+\infty }\frac{[a_k^{(m)}]^j}{j}\sin ^{2j}(u/2) \nonumber
 \\& = - \sum_{j=1}^{+\infty }\frac{[a_k^{(m)}]^j}{j4^j} \left\{\binom{2j}{j} + 2\sum_{s=1}^{j}\left( -1\right)^{s}\binom{2j}{j-s} \cos(su)\right\}
\end{align}
where the second equality follows from \eqref{Line}. Using the duplication formula
\begin{equation*}
(2j)! = 4^j(1/2)_jj!,
\end{equation*}
we get
\begin{align}\label{E2}
\sum_{j=1}^{+\infty }\frac{[a_k^{(m)}]^j}{j4^j} \binom{2j}{j} &= \sum_{j \geq 1}\frac{(1/2)_j}{j!} \frac{[a_k^{(m)}]^j}{j} = \int_0^{a_k^{(m)}}\sum_{j \geq 1}\frac{(1/2)_j}{j!} x^j\frac{dx}{x} \nonumber
\\& = \int_0^{a_k^{(m)}} \frac{1}{\sqrt{1-x}(\sqrt{1-x}+1)}dx = 2\ln\left(\frac{2}{\sqrt{1-a_k^{(m)}}+1}\right).
\end{align}
Now consider the series
\begin{equation*}
\sum_{j \geq s} \binom{2j}{j-s} \frac{[a_k^{(m)}]^j}{j4^j} = \int_0^{a_k^{(m)}}\sum_{j \geq s} \binom{2j}{j-s} \frac{x^j}{4^j} \frac{dx}{x}
\end{equation*}
for fixed $s \geq 1$. Performing an index change $j \mapsto j+s$ and using again the duplication formula
\begin{equation*}
(2j+2s)! = \frac{4^{j+s}}{\sqrt{\pi}} \Gamma(j+s+1/2) (j+s)!,
\end{equation*}
we get
\begin{align*}
\sum_{j \geq s} \binom{2j}{j-s} \frac{x^j}{4^j} &= \frac{\Gamma(s+1/2)\Gamma(s+1)}{\sqrt{\pi}\Gamma(2s+1)} \sum_{j\geq 0}\frac{(s+1/2)_j(s+1)_j}{(2s+1)_j j!}x^{j+s}
\\& = \frac{x^s}{4^s} {}_2F_1(s+1/2, s+1, 2s+1; x)
\end{align*}
where ${}_2F_1$ is the Gauss hypergeometric function (\cite{Erd}, CH.II).
Set
\begin{equation*}
\alpha(x) := \frac{x}{(1+\sqrt{1-x})^{2}} = \frac{1-\sqrt{1-x}}{1+\sqrt{1-x}}, \quad x \in [0,1],
\end{equation*}
recall formula (6), p.101 in \cite{Erd}
\begin{equation*}
{}_2F_1(s+1/2, s+1, 2s+1; x) = \frac{4^s}{\sqrt{1-x}(1+\sqrt{1-x})^{2s}}.
\end{equation*}
Then
\begin{align*}
\sum_{j \geq s} \binom{2j}{j-s} \frac{x^j}{4^j} &= \frac{\alpha(x)^s}{\sqrt{1-x}}.
\end{align*}
But
\begin{equation*}
\alpha'(x) = \frac{\alpha(x)}{x\sqrt{1-x}},
\end{equation*}
which implies
\begin{equation*}
\sum_{j \geq s} \binom{2j}{j-s} \frac{[a_k^{(m)}]^j}{j4^j} = \frac{1}{s} [\alpha^{s}(x)]_0^{a_k^{(m)}} = \frac{1}{s}\alpha^s(a_k^{(m)}).
\end{equation*}
As a result
\begin{equation}\label{E3}
2\sum_{j\geq 1}\frac{[a_k^{(m)}]^j}{j4^j} \sum_{s=1}^{j}\left( -1\right)^{s}\binom{2j}{j-s} \cos(su) = 2\sum_{s \geq 1}\frac{(-1)^s}{s}\alpha^s(a_k^{(m)}) \cos(su).
\end{equation}
The RHS of \eqref{E3} can be written as
\begin{equation*}
\int e^{iux} \left\{\sum_{s \geq 1}\frac{(-1)^s}{s}\alpha^s(a_k^{(m)})[\delta_s+\delta_{-s}](x)\right\} := -\int e^{iux} \mu_k^{(m)}(dx).
\end{equation*}
The signed measure $\mu_k^{(m)}$ is finite:
\begin{equation*}
\int \mu_k^{(m)}(dx) = 2\sum_{s \geq 1}\frac{(-1)^{s-1}}{s}\alpha^s(a_k^{(m)}) = \ln[1+\alpha(a_k^{(m)})] = 2\ln\left(\frac{2}{1+\sqrt{1-a_k^{(m)}}}\right),
\end{equation*}
has finite moments of all orders and obviously satisfies
\begin{equation*}
\int x{\bf 1}_{|x| <1} \mu_k^{(m)}(dx) = 0.
\end{equation*}
Combining \eqref{E1}, \eqref{E2} and \eqref{E3}, the proposition is proved.
\end{proof}

\section{An infinitely divisible distribution}

The preceding computations show that
\begin{equation*}
\ln \mathbb{E}[e^{iuX_m}] = imu + \lambda \int (e^{iux} - 1)\left\{\delta_1
+ \sum_{k=1}^m \mu_k^{(m)}\right\}(dx)
\end{equation*}
which implies again that $X_m$ is not infinitely divisible since the
representative measure is signed. Nonetheless, we have seen that apart from
positivity, $\mu_k^{(m)}$ share the same properties \eqref{LM} of L\'evy
measures. Since its total variation
\begin{equation*}
|\mu_k^{(m)}| := \sum_{s \geq 1}\frac{1}{s}\alpha^s(a_k^{(m)})[\delta_s+%
\delta_{-s}](x)
\end{equation*}
does so and is positive, then
\begin{equation*}
\delta_1 + \sum_{k=1}^m|\mu_k^{(m)}|
\end{equation*}
is a L\'evy measure and the corresponding L\'evy-Kintchine representation
\begin{equation*}
imu + \lambda \int (e^{iux} - 1 - iux\mathbf{1}_{|x| < 1})\left\{\delta_1 +
\sum_{k=1}^m |\mu_k^{(m)}|\right\}(dx)
\end{equation*}
gives rise to an infinitely divisible probability distribution. Using the
linearization formula
\begin{equation}  \label{Linea}
\cos ^{2j}(u/2)=\frac{1}{2^{2j}}\left\{\binom{2j}{j} + 2\sum_{s=1}^{j}\binom{%
2j}{j-s} \cos(su)\right\}.
\end{equation}
( see \cite{Gra-Ryz}, p.31) and reading backward the proof of the previous
proposition, we deduce that the characteristic function of this new
distribution reads
\begin{equation*}
\exp \left( \lambda \left(e^{iu}-1\right)\right) e^{imu} \frac{1}{%
L_{m}^{\left( 0\right) }\left(2\lambda (1+\cos u) \right)}.
\end{equation*}
Recalling $a_k^{(m)} = 4\lambda/x_k^{(m)} \in (0,1)$ then the expansion
\begin{align*}
\frac{1}{L_{m}^{\left( 0\right) }\left(2\lambda (1+\cos u) \right)} &=
\prod_{k=1}^m\frac{1}{1-a_k^{(m)}\cos^2u} = \sum_{\tau_1,\dots, \tau_m \geq
0}\prod_{j=1}^m[a_j^{(m)}]^{\tau_j} [\cos u]^{2(\tau_1+\dots+\tau_m)} \\
& = \sum_{\tau_1,\dots, \tau_m \geq 0}\prod_{j=1}^m[(x_j^{(m)}]^{-\tau_j}
[4\lambda \cos u]^{2(\tau_1+\dots+\tau_m)}
\end{align*}
together with the linearization formula \eqref{Linea} show that $u \mapsto %
\left[L_{m}^{\left( 0\right) }\left(2\lambda (1+\cos u) \right)\right]^{-1}$
is the characteristic function of a discrete probability distribution, say $%
\kappa_m$. Besides, if we denote by $\tau = (\tau_1 \geq \tau_2 \geq \dots
\geq \tau_m \geq 0) \in \mathbb{N}^m$ a partition of length $m$ and weight $%
|\tau| := \tau_1+\dots + \tau_m$, then
\begin{align*}
\frac{1}{L_{m}^{\left( 0\right) }\left(2\lambda (1+\cos u) \right)} = m!
\sum_{\tau}\prod_{j=1}^m[(x_j^{(m)}]^{-\tau_j} [4\lambda \cos u]^{2|\tau|}.
\end{align*}
Consequently, the weights of $\kappa_m$ are seen to be
\begin{eqnarray*}
\kappa_m\{0\} = m! \sum_{\tau} \binom{2|\tau|}{|\tau|} \prod_{j=1}^m \left(%
\frac{\lambda}{x_j^{(m)}}\right)^{\tau_j}
\end{eqnarray*}
and for any integer $n \geq 1$,
\begin{eqnarray*}
\kappa_m\{-n\} = \kappa_m\{n\} = 2m! \sum_{\tau, |\tau| \geq n} \binom{%
2|\tau|}{|\tau|-n} \prod_{j=1}^m \left(\frac{\lambda}{x_j^{(m)}}%
\right)^{\tau_j}.
\end{eqnarray*}

{\small Acknowledgments. }{\footnotesize A major part of this work was
completed when Z. Mouayn was visiting IHES and IRMAR in 2014. He would like
to thank both institutes for their supports and hospitality.}

\end{document}